\documentclass[11pt]{article}
\newtheorem{theorem}{\indent Theorem}[section]
\newtheorem{proposition}[theorem]{\indent Proposition}

\newtheorem{definition}[theorem]{\indent Definition}
\newtheorem{corollary}[theorem]{\indent Corollary}
\newtheorem{example}[theorem]{\indent Example}
\def \N{I\!\!N}

\def \Z{Z\!\!\!Z}
\def \FF{I\!\!F}
\def \P{I\!\!P}
\def \one{{\bf 1}}
\newcommand{\ee}{e^{\frac1e}}

\begin{document}

\title{Stochastic chains with memory of variable length}

\author{A.~Galves\thanks{
This work is part of PRONEX/FAPESP's project
  \emph{Stochastic behavior, critical phenomena and rhythmic pattern
    identification in natural languages} (grant number 03/09930-9),
  CNPq's project \emph{Stochastic modeling of speech} (grant number
  475177/2004-5) and CNRS-FAPESP project \emph{Probabilistic phonology
    of rhythm}. AG is partially supported by a CNPq fellowship (grant
  308656/2005-9). 
  } \and E.~L\"ocherbach }

\maketitle
\thispagestyle{empty}
\pagestyle{empty}

\begin{center}
{\it Dedicated to Jorma Rissanen on his 75'th birthday}
\end{center}

\begin{abstract}
  Stochastic chains with memory of variable length constitute an
  interesting family of stochastic chains of infinite order on a
  finite alphabet. The idea is that for each past, only a finite
  suffix of the past, called {\sl context}, is enough to predict the
  next symbol.  These models were first introduced in the information
  theory literature by Rissanen (1983) as a universal tool to perform
  data compression. Recently, they have been used to model up
  scientific data in areas as different as biology, linguistics and
  music.  This paper presents a personal introductory guide to this
  class of models focusing on the algorithm Context and its rate of
  convergence.
\end{abstract}

\section{Introduction}

Chains with memory of variable length appear in Rissanen's 1983 paper
called {\sl A universal data compression system}. His idea was to
model a string of symbols as a realization of a stochastic chain where
the length of the memory needed to predict the next symbol is not
fixed, but is a deterministic function of the string of the past
symbols.

Considering a memory of variable length is a practical way to overcome
the well known difficulty of the exponentially growing number of
parameters which are needed to describe a Markov chain when its order
increases.  However if one wants to fit accurately complex data using
a Markov chain of fixed order, one has to use a very high order. And
this means that to estimate the parameters of the model we need huge
samples, which makes this approach unsuitable for many practical
issues.

It turns out that in many important scientific data, the length of the
relevant portion of the past is not fixed, on the contrary it depends
on the past. For instance, in molecular biology, the translation of a
gene into a protein is initiated by a fixed specific sequence of
nucleotide bases called {\it start codon}.  In other words, the start
codon designs the end of the relevant portion of the past to be
considered in the translation.

The same phenomenon appears in other scientific domains. For instance
in linguistics, both in phonology and in syntax, there is the notion
of domains in which the grammar operates to define admissible strings
of forthcoming symbols. In other terms, the boundary of the linguistic
domain defines the relevant part of the past for the processing of the
next linguistic units. Rissanen's ingenious idea was to construct a
stochastic model that generalizes this notion of relevant domain to
any kind of symbolic strings.

To be more precise, Rissanen (1983) called {\it context} the relevant
part of the past. The stochastic model is defined by the set of all
contexts and an associated family of transition probabilities.

Models with memory of variable length are not only less expensive than
the classical fixed order Markov chains, but also much more clever
since they take into account the structural dependencies present in
the data. This is precisely what the set of contexts expresses.

Rissanen has introduced models having memory of variable length as a
universal system of data compression. His goal was to compress in real
time a string of symbols generated by an unknown source. To do this,
we have to estimate at each step the length of the context of the
string observed until that time step, as well as the associated
transition probabilities.

If we knew the contexts, then the estimation of the associated
transition probabilities could be done using a classical procedure
such as maximum likelihood estimation.  Therefore, the main point is
to put hands on the context length. In his seminal 1983 paper,
Rissanen solved this problem by introducing the algorithm {\it
  Context}.  This algorithm estimates in a consistent way both the
length of the context as well as the associated transition
probability.

The class of models with memory of variable length raises interesting
questions from the point of view of statistics. Examples are the rate
of convergence and the fluctuations of the algorithm Context and other
estimators of the model. Another challenging question would be how to
produce a robust version of the algorithm Context.

But also from the point of view of probability theory, this class of
models is interesting. In effect, if the length of the contexts is not
bounded, then chains with memory of variable length are chains of
infinite order. Existence, uniqueness, phase-transitions, perfect
simulation are deep mathematical questions that should be addressed to
in this new and challenging class of models.

Last but not least, models of variable length revealed to be very
performing tools in applied statistics, by achieving in an efficient
way classification tasks in proteomics, genomics, linguistics,
classification of musical styles, and much more.

In what follows we present a personal introductory guide to this class
of models with no attempt to give a complete survey of the subject. We
will mainly focus on the algorithm Context and present some recent
results obtained by our research team.

\section{Probabilistic context trees}

In what follows $A$ will represent a finite alphabet of size $|A|$.
Given two integers $m\leq n$, we will denote by $x_m^n$ the sequence
$(x_m, \ldots, x_n)$ of symbols in $A$. Let $A^*_+$ be the set of all
finite sequences, that is
$$  A^*_+ \,=\, \bigcup_{k=1}^{\infty}\,A^{\{ -k, \ldots , -1 \}}.$$
We shall write $\underbar A  = A^ {\{ \ldots, -n, \ldots , -2, -1 \}},$ and denote by $x_{- \infty }^{-1} $ any element of 
$\underbar A .$

Our main object of interest is what we shall call {\it context length
  function}.
\begin{definition}\label{contextlength}
  A {\it context length function} $l$ is a function $l : A_+^* \to \{ 1, 2, \ldots \} \cup \{ \infty \} $ satisfying the following two properties.\\
  (i) For any $ k \geq 1 ,$ for any $x_{-k}^{-1} \in A_+^* ,$ we have
$$ l(x_{-k}^{-1} ) \in \{ 1, \ldots , k \} \cup \{ + \infty \} .$$
(ii) For any $x_{- \infty }^{-1} \in \underbar A ,$ if $l(x_{-k}^{-1}
) =k $ for some $k \geq 1,$ then
$$
\begin{array}{ll}
 l(x_{-i}^{-1} ) = \infty , &\mbox{ for any $i < k$}\\
 l(x_{-i}^{-1} ) = k, &\mbox{ for any $i > k$} .
\end{array} $$
\end{definition}

Intuitively, given a sequence $x_{- \infty }^{-1} ,$ the function $l$
tells us, at which position in the past we can stop since we have
reached the end of the context. The first condition is a kind of
adaptivity condition. It tells us that we can decide whether the end
of the context has already been reached at step $k$ just by inspecting
the past sequence up to that step. If $l$ equals $+ \infty ,$ we have
to look further back in the past.  The second condition is a
consistency condition. It tells us that once we have reached the bound
of the context, we do not have to look further back in the past, and
that the context of a longer sequence $x_{-i}^{-1},$ $i > k,$ is also
the context of $x_{-k}^{-1} .$ In other terms, once the identification
of the context is made at a given step $k,$ this decision will not be
changed by any further data present in the past before $k.$

By abuse of notation, we shall also call $l$ the natural extension of
the context length function to $\underbar A ,$ given by
$$l( x_{-\infty }^{-1}) = \inf \, \{ k \geq 1 : \, l(x_{-k}^{-1}) <  + \infty   \} ,$$
with the convention that $\inf \emptyset = + \infty .$

\begin{definition}
For any $ x_{-\infty }^{-1} \in \underbar A,$ we shall call $ x_{- l( x_{-\infty }^{-1}) }^{-1}$ the {\it context} associated to $l$ 
of the infinite sequence
$ x_{-\infty }^{-1}.$
\end{definition}

\begin{definition}\label{vlm}
  Let $l$ be a given context length function.  A stationary stochastic
  chain $(X_n)_{n \in \Z }$ taking values in $A$ is a chain having
  memory of variable length, if for any infinite past
  $x_{-\infty}^{-1} \in \underbar A$ and any symbol $a \in A, $ we
  have
\begin{equation}\label{eq:real}
 P \left(X_0 =a | X_{-\infty}^{-1}=x_{-\infty}^{-1}\right)\,= 
 P \left(X_0 =a |
   X_{- l(x_{-\infty}^{-1})}^{-1}=x_{-  l(x_{-\infty}^{-1})}^{-1}\right)\, .\, 
\end{equation}
We shall use the short hand notation  
\begin{equation}\label{petitp}
 p(a | x_{-k }^{-1} ) = P \left(X_0 =a |
   X_{- k}^{-1}=x_{- k }^{-1}\right) .
   \end{equation}
We are mainly interested in those values of $    p(a | x_{-k }^{-1} ) $ 
where $k = l(x_{-\infty}^{-1}) .$
\end{definition}

Observe that
the set  $\{ \ell(X_{-\infty}^{-1}) = k\}$ is measurable with respect
to the $\sigma$-algebra  generated by $X_{-k}^{-1}$. Thus we have 

\begin{proposition}
  Let $(X_n) $ be a stationary chain as in definition \ref{vlm},
  having context length function $l.$ Put $\FF_k = \sigma \{ X_{-k },
  \ldots , X_{-1} \} , k \geq 1 .$ Then $l( X_{-\infty}^{-1})$ is a
  $(\FF_k )_k -$stopping time.
\end{proposition}

Given a context length function $l ,$ we define an associated
countable subset $\tau \subset A^*_+$ by
$$ \tau = \tau^l = \{ x_{-k}^{-1}  \, : k = l( x_{-k}^{-1}) , \, k \geq 1  \} .$$
To simplify notation, we will denote by $\underbar x $ and $\underbar y$ 
generic elements of $\tau .$ 

\begin{definition}
  Given a finite sequence $x_{-k}^{-1},$ we shall call {\it suffix} of
  $x_{-k}^{-1} $ each string $x_{-j}^{-1} $ with $j \le k .$ If $j < k
  $ we call $x_{-j}^{-1} $ {\it proper suffix} of $x_{-k}^{-1} .$ Now
  let ${\mathcal S} \subset A^*_+ .$ We say that ${\mathcal S}$
  satisfies the {\it suffix property} if $\mbox{ for no } x_{-k}^{-1}
  \in {\mathcal S} , \mbox{ there exists a proper suffix } x_{-j}^{-1}
  \in {\mathcal S} \mbox{ of }x_{-k}^{-1} .$
\end{definition}

The following proposition follows immediately from property (ii) of definition \ref{contextlength}. 
\begin{proposition}
Given a context length function $l,$ the associated set $\tau^l $ satisfies the {\it suffix property}. 
\end{proposition}

As a consequence, $\tau $ can be identified with the set of leaves of a rooted tree with a countable set of finite
labeled branches. 

\begin{definition}
We call {\it probabilistic context tree on } $A$ the ordered pair $(\tau , p) ,$ where 
$$ p = \{ p(.| \underbar  x  ) ,  \,\underbar  x  \in \tau \} $$
is the family of transition probabilities of (\ref{petitp}). We say
that the probabilistic context tree $(\tau, p)$ is unbounded if the
function $l$ is unbounded.
\end{definition}

\begin{definition}
  Let $(X_n)_{n \in \Z} $ be a stationary chain and let $(\tau, p)$ be
  a probabilistic context tree.  We shall say that $(X_n)_{n \in \Z }
  $ is {\it compatible} with $(\tau, p) ,$ if (\ref{petitp}) holds for
  all $ \underbar x \in \tau .$
\end{definition}

In order to illustrate these mathematical concepts, let us consider
the following example.

\begin{example}
  Consider a two-symbol alphabet $A = \{ 0, 1 \} $ and the following
  context length function
$$ l(x_{-\infty }^{-1} ) = \inf \, \{ k : x_{-k} = 1 \} .$$
Then the associated tree $\tau $ is given by 
$$ \tau = \{ 1 0^k, k \geq 0 \} ,$$
where $10^k $ represents the sequence $(x_{-k -1}, x_{-k},\ldots ,
x_{-1}) $ such that $x_{-i} = 0 $ for all $1 \le i \le k $ and $x_{-k
  -1} = 1.$

The associated transition probabilities are defined by 
$$ P( X_0 = 1 | X_{-1} = \ldots = X_{-k} =0, X_{-k -1} = 1) = q_k , k \geq 0 ,$$
with $0 \le q_k  \le1 .$ 
\end{example}

Clearly, the stochastic chain associated to this context length
function $l$ is a chain of infinite order. This raises the
mathematical question of existence of such a process. It is
straightforward to see that the following proposition holds true.

\begin{proposition}
  Suppose that $X_0 = 1.$ Put $T_1 = \inf \{ k \geq 1 : X_k = 1 \} .$
  A necessary and sufficient condition for $T_1 < + \infty $ almost
  surely is
\begin{equation}\label{suminfini}
\sum_{k\geq 0}  q_k = + \infty .
\end{equation} 
\end{proposition}
This means that if (\ref{suminfini}) is satisfied, then -- provided
the chain starts from $1$ -- almost surely there will be appearance of
an infinite number of the symbol $1.$ This implies that there exists a
non-trivial stationary chain associated to this probabilistic context
tree.

Observe that the process $X_n$ is actually a renewal process with the
renewal times defined as follows.
$$T_{0} = \sup \{ n< 0 : X_n = 1 \} ,$$
and for $k \geq 1 ,$  
$$T_{k } = \inf \{ n > T_{k-1} : X_n = 1\} \  \mbox{ and }\
 \, T_{-k} := \sup \{ n < T_{-(k-1)} : X_n = 1 \}  . $$ 
 
This example shows clearly that the tree of contexts defines a partition of all possible pasts 
with the exception of the single string composed of all symbols identical to $0.$ The condition (\ref{suminfini})
shows that it is possible to construct the chain taking values in the set of all sequences having an infinite
number of symbols $1$ both to the left and to the right of the origin. 

However, we could also include this
exceptional string to the set of possible contexts by defining an extra parameter $q_\infty .$ This is the
choice of Csisz{\'a}r and Talata (2006). If $q_\infty$ is strictly positive, then 
condition (\ref{suminfini}) implies that after a finite time, there will appearance of the symbol $1,$ even if we start
with an infinity of symbols $0.$  In other terms, this exceptional string does not have to be considered if we
are interested in the stationary regime of the chain.  

In case that $q_\infty = 0 $ and (\ref{suminfini}) holds, we have the phenomenon of phase transition. One of
the phases is composed of only one string having only the symbol $0.$   

The renewal process is an interesting example of a chain having memory of 
unbounded variable length. In the case where the probabilistic context tree is bounded,
the corresponding chain is in fact a Markov chain whose order is equal to the maximal
context length. However, the tree of contexts provides interesting additional
information concerning the dependencies in the data and the structure of
the chain. This raises the issue how to estimate the context tree out of the data. This was
originally solved in Rissanen's 1983 paper using the algorithm Context. 

At this point it is important to discuss the following minimality issue. Among all possible context trees
fitting the data, we want of course to identify the smallest one. This is the tree corresponding
to the smallest context length function. More precisely, if $l$ and $l'$ are context length
functions, we shall say that $l \le l' $ if $l(x_{- \infty }^{-1} ) \le l'(x_{- \infty }^{-1} )$ for any
string $x_{- \infty }^{-1} \in \underbar A .$ From now on we shall call {\it context} of a string
$x_{- \infty }^{-1} $ the context associated to the minimal context length function. Estimating this
minimal context is precisely the goal of the algorithm Context. 

\section{The algorithm Context}
We now present the algorithm Context
introduced by Rissanen (1983). The goal of the algorithm is to estimate adaptively the
context of the next symbol $X_n$ given the past symbols $X_0^{n-1}$. 
The way the algorithm Context works can be
summarized as follows. Given a sample produced by a chain with
variable memory, we start with a maximal tree of candidate contexts
for the sample. The branches of this first tree are then pruned 
starting from the leaves towards the root until
we obtain a minimal tree of contexts well adapted to the sample. We
associate to each context an estimated probability transition defined
as the proportion of time the context appears in the sample followed
by each one of the symbols in the alphabet.  We stop pruning once
the gain function exceeds a given
threshold.

Let $X_0, X_1, \ldots, X_{n-1}$ be a sample from the finite probabilistic
tree $(\tau,p)$.  For any finite string $x_{-j}^{-1}$ with $j \le n$,
we denote $N_n(x_{-j}^{-1})$ the number of occurrences of the string
in the sample
\begin{equation}
  \label{eq:Nn}
  N_n(x_{-j}^{-1}) = \sum_{t=0}^{n-j} \one\left\{X_t^{t+j-1} =
  x_{-j}^{-1}\right\}\,.
\end{equation}
Rissanen first constructs a maximal candidate context $X_{n-M(n)}^{n-1}$ where 
$M(n)$ is a random length defined as follows
\begin{equation}
  \label{eq:mn}
  M(n) \, = \, \min \left\{ i=0,1, \ldots, \lfloor C_1 \, \log n
  \rfloor : N_n (X_{n - i}^{n-1}) \,>\, \frac{C_2 \, n}{\sqrt{\log n}}
  \right\}\, .
\end{equation}
Here $C_1$ and $C_2$ are arbitrary positive constants. In the case 
the set is empty we take $M(n)=0$. 

Rissanen then shortens this maximal candidate context by successively
pruning the branches according to a sequence
of tests based on the likelihood ratio statistics. This is formally
done as follows.

If $\sum_{b \in A}N_n(x_{-k}^{-1}b) \,>\,0$, define 
the estimator of the transition probability $p$  by
\begin{equation}
  \label{eq:phat}
  \hat{p}_n(a|x_{-k}^{-1}) = \frac{N_n(x_{-k}^{-1}a)}{\sum_{b \in
  A}N_n(x_{-k}^{-1}b)} \, 
\end{equation}
where $x_{-j}^{-1}a$ denotes the string $(x_{-j}, \ldots, ,
x_{-1},a)$, obtained by concatenating $x_{-j}^{-1}$ and the symbol
$a$. If $\sum_{b \in A}N_n(x_{-k}^{-1}b) \,=\,0$, define
 $\hat{p}_n(a|x_{-k}^{-1}) \,=\, 1/|A|$.  

For  $i  \geq 1$ we define
 \begin{equation}
  \label{eq:delta}
\Lambda_n (x_{- i }^{-1} )
\,=\,  2 \, \sum_{y \in A} \sum_{a \in A}  N_n(y x_{-i}^{-1}a)
\log\left[\frac{\hat{p}_n(a|x_{-i}^{-1}y)}
{\hat{p}_n(a|x_{-i}^{-1})} \right]\, ,
\end{equation}
where $y x_{-i}^{-1} $ denotes the string $(y,x_{-i}, \ldots ,  x_{-1}), $ 
and where 
$$\hat{p}_n(a|x_{-i}^{-1}y)= \frac{N_n(y x_{-i}^{-1}a)}{\sum_{b \in
  A}N_n(y x_{-i}^{-1}b)}. $$

Notice that $\Lambda_n (x_{- i }^{-1})$ is the log-likelihood ratio
statistic for testing the consistency of the sample with a
probabilistic suffix tree $(\tau,p)$ against the alternative that it is
consistent with $(\tau',p')$ where $\tau$ and $\tau'$ differ only by
one set of sibling nodes branching from $x_{-i}^{-1}$. 
$\Lambda_n (x_{- i }^{-1})$ plays the role of a gain function telling
us whether it is worth or not taking a next step further back in the past.

Rissanen then defines the length of the estimated current context $\hat{\ell}_n$
as
\begin{equation}
\label{eq:ell}
\hat{\ell}_n(X_0^{n-1})= 1+ \max \left\{i=1,\ldots, M(n)-1: \Lambda_n
(X_{n-i}^{n-1}) \,>\, C_2 \log n \right\}\, ,
\end{equation}
where $C_2$ is any positive constant.
\vspace{.5cm}

Then, the result in Rissanen (1983) is the following.

\begin{theorem}\label{teo0}
Given a realization $X_0,
\ldots, X_{n-1}$ of a probabilistic suffix tree $(\tau,p)$ with finite height,
then 
\begin{equation}
  \label{eq:riss}
P\left( \hat{\ell}_n(X_0^{n-1}) \neq \ell(X_0^{n-1}) \right) 
\longrightarrow 0  
\end{equation}
as $n \rightarrow \infty$. 
\end{theorem}

Rissanen proves this result in a very short and elegant way. 
His starting point is the following upper bound.
\begin{eqnarray}
  \label{eq:riss2}
&&\lefteqn{P\left( \hat{\ell}_n(X_0^{n-1}) \neq \ell(X_0^{n-1}) \right)
  \le} \nonumber \\
&&  P\left(\hat{\ell}_n(X_0^{n-1}) \neq
  \ell(X_0^{n-1})  |   N_n \! \! \left(X_{n-\ell(X_0^{n-1})}^{n-1} \right) >
  \frac{C_2 n}{\sqrt{\log n}} \right) \! P\left(\! \! N_n
  \left(X_{n-\ell(X_0^{n-1})}^{n-1} \right) > \frac{C_2 n}{\sqrt{\log n}}
  \right) \nonumber \\
&&\quad \quad \quad \quad + P \left(\bigcup_{w \in \tau} \left\{
    N_n \left(w \right) \le \frac{C_2 n}{\sqrt{\log n}} \right\}
\right) \, .
\end{eqnarray}

Then he provides the following explicit
upper bound for the conditional probability in the right-hand side of
(\ref{eq:riss2})
\begin{equation}
  \label{eq:riss_bd}
  P\left( \hat{\ell}_n(X_0^{n-1}) \neq
  \ell(X_0^{n-1}) | N_n \left(X_{n-\ell(X_0^{n-1})}^{n-1} \right) >
  \frac{C_2 n}{\sqrt{\log n}} \right) \le  C_1 \, \log n
\, e^{-C'_2 \sqrt{\log n}}  \, ,
\end{equation}
where $C_1$, $C_2$ and $C'_2$ are positive constants independent of the
maximum of the context length function.  

With respect to the second term he only observes that, by ergodicity,
for each $x_{-k}^{-1}\in \tau$ we have
\begin{equation}
  \label{eq:riss_2}
  P\left( N_n \left(x_{-k}^{-1} \right) \le \frac{C_2
  n}{\sqrt{\log n}} \right) \longrightarrow 0
\end{equation}
as $n \rightarrow \infty$. Since $\tau$ is finite the convergence in
(\ref{eq:riss_2}) implies the desired result.

\section{The unbounded case}
In his original paper, Rissanen was only interested in the case of 
bounded context trees. 
However, from the mathematical point of view, it is interesting to
consider also the case of unbounded probabilistic context trees 
corresponding to chains of infinite order. It can be argued that
also from an applied point of view the unbounded case must be 
considered as noisy observation of Markov chains generically 
have infinite order memory. 

The unbounded case raises immediately the preliminary question 
of existence and uniqueness of the corresponding chain. 
This issue can be addressed by adapting to probabilistic context
trees the conditions for existence and uniqueness that have 
already been proved for infinite order chains. This is precisely 
what is done in the paper by Duarte et al. (2006) who adapt
the type A condition presented in Fern\'andez and Galves (2002)
in the following way. 

To simplify the presentation, let us introduce some extra notation. 
Recall that  $\underbar x $ and $\underbar y$ 
denote generic elements of $\tau .$ Given $\underbar x = x_{-i}^{-1} $ and
$\underbar y = y_{-j}^{-1} $, we shall write $\underbar x \stackrel{k}{=}
\underbar y $ if and only if $k \le \min \{ i, j\} $ and $x_{-1} = y_{-1},
\ldots , x_{-k } = y_{-k} .$   

\begin{definition} 
A probabilistic suffix tree $(\tau,p)$ on $A$ is of {\it type
    A} if its transition probabilities $p$ satisfy the following
  conditions.
\begin{enumerate}
\item  {\bf Weakly non-nullness}, that is 
  \begin{equation}
    \label{eq:weak}
    \sum_{a \in A} \inf_{{\tiny \underbar x}  \in \tau} p(a | \underbar x) > 0\, ; 
  \end{equation}
\item {\bf Continuity}, that is
  \begin{equation}
\label{eq:cont}
\beta(k) \,  =  \,  \max_{a \in A} \sup \{ |p(a | \underbar x)
- p(a | \underbar y)|, \underbar y \in \tau, \underbar x \in \tau \, \mbox{ with } \,
\underbar x \stackrel{k}{=} \underbar  y   \} \, \rightarrow  \, 0
\end{equation}  
as $k \rightarrow \infty$. We also define 
$$\beta(0) =  \max_{a \in A} \sup \{ |p(a | \underbar x)
- p(a | \underbar y)|, \underbar y \in \tau, \underbar x \in \tau \, \mbox{ with } \,
x_{-1} \neq y_{-1}  \}.$$

The sequence $\{\beta(k)\}_k \in \N$ is called
the {\bf continuity rate}. 
\end{enumerate}
\end{definition}

For a probabilistic suffix tree of type A with summable continuity
rate, the maximal coupling argument used in Fern\'andez and Galves (2002)
implies the uniqueness of the law of the chain consistent with it.

We now present a slightly different version 
of the algorithm Context using the same gain function $\Lambda_n$  
but in which the length of the maximum context candidate
is now deterministic and nor more random. More precisely, 
we define the length of the biggest candidate context now as 
\begin{equation}\label{detlength}
k(n)= C_1 \log n
\end{equation}
with a suitable positive constant $C_1$. 

The intuitive reason behind the
choice of the upper bound length $ C_1 \, \log n$ is the impossibility
of estimating the probability of sequences of length much longer than
$\log n$ based on a sample of length $n$. Recent versions of this fact
can be found in Marton and Shields (1994, 1996)
and  Csisz{\'a}r (2002). 

Now, the definition of $\hat{\ell}_n$ is similar to the one 
in the original algorithm of Rissanen, that is
\begin{equation}
\label{eq:ell3}
\hat{\ell}_n(X_0^{n-1})= 1 +\max \left\{i=1,\ldots, k(n)-1: \Lambda_n
(X_{n-i }^{n-1}) \,>\, C_2 \log n \right\}\, ,
\end{equation}
where $C_2$ is any positive constant.
\vspace{.5cm} 

The reason for taking the length of the maximum context candidate deterministic
and no more random is to be able to use the classical results on the convergence of
the law of $\Lambda_n(x_{-i}^{-1} )$ to a chi-square distribution. 
However, we are not in a Markov setup since the probabilistic context tree
is unbounded, and the chi-square approximation only works for Markov chains of fixed finite
order. 

To overcome this difficulty, we use the canonical Markov approximation 
of chains of infinite order presented in Fern\'andez and Galves (2002) that we
recall now by adapting the definitions and theorem to the framework of
probabilistic context trees. The goal is to approximate a chain compatible
with an unbounded probabilistic context tree by a sequence of chains 
compatible with bounded probabilistic context trees.

\begin{definition} 
For all $k \ge 1 ,$ the {\it canonical Markov approximation of order $k$} of
  a chain $(X_n)_{n \in \Z}$ is the chain with memory of
  variable length bounded by $k$ compatible with the 
  probabilistic context tree 
  $(\tau^{[k]},p^{[k]})$ where 
\begin{equation}
  \label{eq:cano}
\tau^{[k]} = \{\underbar x \in \tau; l(\underbar x) \le k\} \cup \{ x_{-k}^{-1}; \underbar x \in
\tau, l(\underbar x )\ge k\}
\end{equation}
for all $a \in {A},$ $\underbar x \in \tau,$  and where 
\begin{equation}
    \label{eq:201}
    p^{[k]}(a | x_{-j}^{-1}) \,:=\, P (X_0 = a | X_{-j}^{-1}
    = x_{-j}^{-1}) 
  \end{equation}
  for all $x_{-j}^{-1} \in \tau^{[k]} .$
  \end{definition}

Observe that for contexts $\underbar x \in \tau$ which length does not
exceed $k$, we have $p^{[k]}(a | \underbar x) = p(a | \underbar x)$. However,
for sequences $x_{-k}^{-1}$ which are internal nodes of $\tau$, there
is no easy explicit formula expressing $p^{[k]}(\cdot |x_{-k}^{-1})$
in terms of the family $\{p(\cdot| \underbar y), \underbar y  \in \tau\}$. \\

The main result of Fern\'andez and Galves (2002) that will be used in the proof 
of the consistency of the algorithm Context  can be stated as follows.

\begin{theorem}
 Let $(X_n)_{n \in \Z}$ be a chain compatible
with a type A probabilistic context tree $(\tau,p)$ with summable
continuity rate, and let $(X^{[k]}_n)_{n \in \Z}$ be its canonical Markov
approximation of order $k$. Then there exists a coupling between
$(X_n)_{n\in \Z}$ and $(X^{[k]}_n)_{n \in \Z}$ and a constant $C > 0$ such that
  \begin{equation}
    \label{eq:FG6}
    P \left( X_0 \neq X^{[k]}_0 \right)
    \le C  \beta(k)\, . 
  \end{equation}
\end{theorem}
\vspace{.5cm}

Using this result and the classical chi-square approximation for Markov
chains, Duarte et al. (2006) proved the consistency of their 
version of the algorithm Context in the unbounded case and also
provided an upper bound for the rate of convergence. Their result is the
following. 

\begin{theorem} \label{teo1} Let  $X_0, X_2, \ldots, X_{n-1}$ be a sample from
  a type A unbounded probabilistic suffix tree $(\tau,p)$ with
  continuity rate $\beta(j)\le f(j) \exp\{-j\}$, with $f(j)
  \rightarrow 0$ as $j \rightarrow \infty$. Then, for any choice of
  positive constants $C_1$ and $C_2$ in (\ref{detlength}) and (\ref{eq:ell3}),
  there exist positive constants $C$ and $D$ such that
$$
P\left(\hat{\ell}_n(X_0^{n-1}) \neq \ell(X_0^{n-1}) \right) \le
C_1 \log n (n^{-C_2 } + D/n)  + C f(C_1 \log n) \, .$$
\end{theorem}

The proof can be sketched very easily. Take $k=k(n)=C_1 \log(n)$
and construct a coupled version of the processes $(X_t)_{t \in \Z} $
and $(X^{[k(n)]}_t)_{t \in \Z }.$ First of all notice that for $k = k(n),$ 
\begin{eqnarray}\label{upperbound1}
&& P\left(\hat{\ell}_n(X_0, \ldots, X_{n-1}) \neq  \ell(X_0, \ldots , X_{n-1}) \right) \le \nonumber \\
&&
P \left(\hat{\ell}_n(X_0^{[k]}, \ldots X_{n-1}^{[k]} ) \neq \ell(X_0^{[k]}, \ldots X_{n-1}^{[k]} ) \right)  + P \left(
  \bigcup_{i=1}^{n} \{X_i \neq X_i^{[k]}\} \right). 
  \end{eqnarray}
Using the inequality (\ref{eq:FG6}) of Fern\'andez and Galves (2002), the second
term in (\ref{upperbound1}) can be bounded above as 
$$ P \left(  \bigcup_{i=1}^{n} \{X_i \neq X_i^{[k]}\} \right)  \le n \, C \, \beta(k(n)) .  $$ 

The first term in (\ref{upperbound1}) can be treated using the 
classical chi-square approximation for the log-likelihood ratio
test for Markov chains of fixed order $k.$ 

More precisely, we 
know that for fixed $x_{-i}^{-1}$, under the null hypothesis, the
statistics $\Lambda_n (x_{-i}^{-1})$, given by (\ref{eq:delta}), has
asymptotically chi-square distribution with $|A| -1$ degrees of
freedom (see, for example, van der Vaart (1998)). We recall that,
for each $x_{-i}^{-1}$ the null hypothesis ($H^i_0$) is that the true context
is $x_{-i}^{-1}$.

Since we are going to perform a sequence of $k(n)$ sequential tests
where $k(n) \rightarrow \infty$ as $n$ diverges, we need to control
the error in the chi-square approximation. For this, we use a
well-known asymptotic expansion for the distribution of $\Lambda_n
(x_{-i}^{-1})$ due to Hayakawa (1977) which implies that
\begin{equation}
  \label{eq:haya}
 P \left ( \Lambda_n (x_{-i}^{-1}) \le t  | H^i_0 \right) = P \left(
  \chi^2 \le t \right) + D/n \, , 
\end{equation}
where $D$ is a positive constant and $\chi^2$ is random variable with
distribution chi-square with $|A| -1$ degrees of freedom.

Therefore, it is immediate that
$$ P \left(\Lambda_n (x_{-i}^{-1}) \,>\,  C_2 \log n \right) \le e^{-C_2 \log
  n} + D/n \,.$$

By the way we defined $\hat{\ell}_n$ in (\ref{eq:ell3}), in order to find $\hat{\ell}_n(X_0^{n-1})$ we have to
perform at most $k(n)$ tests. We want to give an upper bound for the
overall probability of type I error in a sequence of $k(n)$ sequential
tests. An upper bound is given by the Bonferroni inequality, which
in our case can be written as
$$
P \left( \cup_{i=2}^{k(n)} \{ \Lambda_n (x_{-i}^{-1}) > C_2 \log n \}|
  H^i_0 \right) \le \sum_{i=2}^{k(n)} P ( \Lambda_n (x_{-i}^{-1}) > C_2
  \log n | H^i_0 ).$$
This last term is bounded above by $ C_1 \log n (n^{- C_2} +
  D/n)$. This concludes the proof. 
  
Theorem \ref{teo1} not only proves the consistency of the algorithm Context,
but it also gives an upper bound for the rate of convergence. 
The estimation of the rate of convergence is crucial because it gives a bound on the
minimum size of a sample required to guarantee, with a given probability, that the
estimated tree is the good one. This is the issue we
address to in the next section. 

\section{Rate of convergence of the algorithm Context}
Note that Rissanen's original theorem \ref{teo0} as well as
theorem \ref{teo1} only show that
all the contexts identified are true contexts with high probability.
In other words, the estimated tree is a subtree of the true
tree with high probability. In the
case of bounded probabilistic context trees
this missing point was handled with in Weinberger et al. (1995). 
This paper not only proves that the set of all contexts is 
reached, but also gives a bound for the rate of convergence. 

More precisely, let us define the empirical tree 
\begin{equation}\label{emptree}
\hat{\tau}_n = \left\{ X^{j-1}_{j- \hat{\ell}_j (X_{0}^{j-1})} : j = n/2 , \ldots , n  \right\} .
\end{equation}
Actually, this is a slightly simplified version of the empirical tree defined in 
Weinberger et al. (1995). In particular, we are neglecting all the computational 
aspects considered there. But from the mathematical point of view, this 
definition perfectly does the job. Their convergence result is the following.

\begin{theorem}\label{weinbergerlemma1}
Let $(\tau , p)$ be a bounded probabilistic context tree and let $X_0, \ldots , X_n$
be compatible with $(\tau , p ) .$ Then we have 
$$ \sum_{n\geq 1 } P( \hat{\tau}_n \neq \tau ) \log n < + \infty .$$
\end{theorem} 

In the unbounded case, this issue was treated without estimation of the rate 
of convergence in Ferrari and Wyner (2003) and including estimation of the
rate of convergence in Galves and Leonardi (2008). 
 
This last paper considers another slightly modified version of the algorithm Context
using a different gain function, 
which has been introduced in Galves et al. (2007). More precisely,  
let us define for any finite string $x_{-k}^{-1} \in A_+^*$ the gain function
\[ 
\Delta_n(x_{-k}^{-1}) = \max_{a\in A}
|\hat{p}_n(a|x_{-k}^{-1})-\hat{p}_n(a| x_{-(k-1)}^{-1})| \/.
\] 
This gain function is well adapted to use exponential inequalities for 
the empirical transition probabilities in the pruning procedure rather
than the chi-square approximation of the log-likelihood ratio as in 
theorems \ref{teo0} and \ref{teo1}. 

The theorem is stated in the following framework. 
Consider a  a stationary chain $(X_n)_{ n \in \Z}$ compatible with
an unbounded probabilistic context tree $(\tau , p).$ For this chain,
we define the sequence $(\alpha_n)_{n \in \N } $ by 
\begin{eqnarray*}
\alpha_0 &=& \sum_{a \in A } \inf_{ { \tiny \underbar x} \in \tau } p(a | \underbar x ) ,  \\
\alpha_{n} &=& \inf_{x_{-n}^{-1} 
} \; \sum_{a \in A } \;  \; \; \inf_{{\tiny \underbar y}  \in \tau : l({\tiny \underbar y}) \geq n, {\tiny \underbar y} \stackrel{n}{=} x_{-n}^{-1}  } \; \; p(a| \underbar y ) .
\end{eqnarray*}
We assume that the probabilistic context tree $(\tau , p)$
satisfies the condition (\ref{eq:weak}) of weakly non-nullness, that is 
$\alpha_0 > 0 .$ We assume also the following summability condition
\begin{equation}\label{eq:alpha}
\alpha = \sum_{n \geq 0 } (1 - \alpha_n)  < +\infty .
\end{equation}

Given a sequence $x_1^j = (x_1, \ldots , x_j)  \in A^j$ we denote by 
$$ p(x_1^j) \,=\,  \P(X_1^j = x_1^j) .$$

Then for an integer $m \geq 1, $ we define
\begin{equation}\label{dk}
D_m = \min_{x_{-k}^{-1} \in\tau: k \leq m}\; \max_{a\in A}\{\,|p(a| x_{-k}^{-1}) - p(a| x_{-(k-1)}^{-1})|\,\},
\end{equation}
 and 
\begin{equation}\label{ek}
\epsilon_m= \min\{\,p(x_{-k}^{-1} )\colon k \leq m \mbox{ and } p(x_{-k}^{-1}) > 0\,\}.
\end{equation}
Intuitively, $D_m$ tells us how distinctive is the difference between 
transition probabilities associated to the exact contexts and those
associated to a string shorter one step in the past. We do not want to
impose restrictions on the transition probabilities elsewhere then at
the end of the branches of the context tree. This has to do with the pruning
procedure which goes from the leaves to the root of the tree. 

In the unbounded case, a natural way to state the convergence results 
is to consider truncated trees. The definition is the following. 
Given an integer $K$ we will denote by $\tau|_K$ the
tree $\tau$ \emph{truncated} to level $K$, that is
$$  \tau|_K = \{ x_{-k}^{-1}  \in \tau\colon k  \le K\} \cup \{ x_{-K}^{-1}  
  \mbox{ such that  } x_{k}^{-1} \in \tau \mbox{ for some } k \geq K \}.$$
Actually, this is exactly the same tree which was called $\tau^{[K]}$ in (\ref{eq:cano}).
The notation $  \tau|_K  $ is more suitable for what follows.

The associated empirical tree of height $k$ is defined in the following way. 
\begin{definition}\label{estim-tree}
  Given $\delta > 0$ and $k < n$, the empirical tree is defined as 
  \[
  \hat{\tau}_n^{\delta, k} = \{x_{-r}^{-1} , 1 \le r \le k :  \Delta_n(x_{-r}^{-1} ) >
  \delta\, \wedge \,\Delta_n( y_{-(r+j)}^{-(r+1)} x_{-r}^{-1} )\leq \delta  , \; \forall \;  
  y_{-j}^{-(r+1)} , 1 \le j \le k-r \} . 
  \]
 In case $r = k ,$ the string $  y_{-(r+j)}^{-(r+1)} $ is empty. 
 \end{definition}
Note that in this definition, the parameter $\delta $ expresses the coarseness of 
the pruning criterion and $k$ is the maximal length of the estimated contexts.

Now, Galves and Leonardi (2008) obtain the following result on the rate
of convergence for the truncated context tree.

\begin{theorem}\label{expobounded}
 Let $(\tau,p)$ be a probabilistic context tree satisfying (\ref{eq:weak})
 and (\ref{eq:alpha}).  Let $X_0, \ldots , X_n $ be a stationary 
 stochastic chain compatible with $(\tau , p).$ 
 Then for any integer $K$, any $k$ satisfying 
 \begin{equation}\label{d}
 k >  \max_{{\tiny \underbar x }\in \tau|_K} \min\,\{ \ell(\underbar y)\colon \underbar y \in\tau,  \; \underbar x \stackrel{K}{=} \underbar y\},
 \end{equation}
 for any 
 $\delta < D_k$ and for  each
 \begin{equation}\label{eq:n}
 n >  \frac{2(|A|+1)}{\min(\delta,D_k-\delta)\epsilon_k}+k
 \end{equation}
 we
  have that 
 $$ P(\hat\tau_n^{\delta,k}|_K \neq \tau|_K)\;\; \leq 
    4\,\ee \,  |A|^{k+2} \exp \bigl[- (n\!-\!k\!)\;
\frac{[\min(\frac{\delta}{2},\frac{ D_k-\delta}{2})-\frac{|A|+1}{(n\!-k\!)\epsilon_k}]^2 
\epsilon_k^2 C}{4 |A|^2(k+1)}\bigl],$$
where
  \[
  C = \frac{\alpha_0}{8e(\alpha+\alpha_0 )}.
  \]
\end{theorem}
In this theorem, the empirical trees have to be of height $k \geq K $ for the 
following reason. Truncating $\tau $ at level $K$ implies that contexts longer than $K$ are 
cut before reaching their end, and associated transition probabilities might not differ
when comparing them at length $K$ and $K-1 .$ That's why we consider the bigger 
empirical tree of height $k$ satisfying condition (\ref{d}). This guaranties that for each
element $\underbar x $ of the truncated empirical tree there is at least one real context $\underbar y$
which has $\underbar x $ as its suffix.
   
As a consequence of theorem \ref{expobounded}, Galves and Leonardi (2008)
obtain the following strong consistency result.

\begin{corollary}\label{main_cor}
Let $(\tau,p)$ be a probabilistic context tree satisfying the conditions of theorem \ref{expobounded}.
Then
  \begin{equation}
  \hat\tau_n^{\delta,k}|_K = \tau|_K,
  \end{equation}
  eventually almost surely, as $n \to \infty .$
\end{corollary}   
   
The main ingredient of the proof of theorem \ref{expobounded}
is an exponential upper bound for the deviations of the 
empirical transition probabilities. More precisely, Galves and Leonardi (2008)
prove the following result. 

\begin{theorem}\label{cor:estim1} 
 For any finite sequence $x_{-k}^{-1}$ with $p(x_{-k}^{-1})>0$, any symbol $a\in A$,
  any $t>0$ and any $n > \frac{|A|+1}{tp(x_{-k}^{-1})}+k$ the following inequality holds.
\begin{eqnarray}\label{expo2}
 & &P\bigl(|\hat{p}_n(a|x_{-k}^{-1})-p(a|x_{-k}^{-1})|>t \bigl)\;\leq  \nonumber\\
 & & \quad \quad \quad \quad 2|A| \, \ee \exp \bigl[- (n\!-\!k\!)\;\frac{[t-\frac{|A|+1}{(n\!-\!k\!)p(x_{-k}^{-1})}]^2 p(x_{-k}^{-1})^2C}{4|A|^2(k+1)}\bigl],
\end{eqnarray}
where 
\begin{equation}\label{C}
    C = \frac{\alpha_0}{8e(\alpha +\alpha_0 )}.
  \end{equation}
\end{theorem}

The proof of this theorem is inspired by recent exponential upper bounds obtained by
Dedecker and Doukhan (2003),
Dedecker and Prieur (2005) and Maume-Deschamps (2006). It is based on the following 
loss-of-memory inequality of Comets et al. (2002). 

\begin{theorem}\label{mixing}
Let $(X_n)_{ n \in\Z }$ be a stationary stochastic chain compatible with the
  probabilistic context tree $(\tau,p)$ of theorem \ref{expobounded}. Then, 
  there exist a sequence $\{\rho_l\}_{l\in\N}$ such that
  for any $i\geq 1$, any $k > i$, any $j\geq 1$ and any
  finite sequence $x_1^j$, the following inequality holds
  \begin{equation}\label{eqmixing}
    \sup_{x_1^{i}\in A^i} 
    | P(X_{k}^{k+j-1}=x_1^j | X_1^i=y_1^i)-p(x_1^j)|\;
    \leq \;  j \, \rho_{k-i-1}\,.
  \end{equation}
  Moreover, the sequence $\{\rho_l\}_{l\in\N}$ is summable and
  \[
  \sum_{l\in\N} \rho_l \;\leq\; 1+\frac{2 \alpha}{\alpha_0}.
  \]
\end{theorem}

Theorem \ref{expobounded} generalizes to the unbounded case previous
results in Galves et al. (2008) for the case of bounded context trees.
Note that the definition of the context tree estimator depends on the
parameter $\delta$, the same appearing in the constants of the
exponential bound.  To assure the consistency of the estimator we have
to choose a $\delta$ sufficiently small, depending on the true
probabilities of the process.  The same thing happens to the parameter
$k.$ Therefore, this estimator is not universal, meaning that for
fixed $\delta$ and $k$ it fails to be consistent for all variable
memory processes for which conditions (\ref{d}) and (\ref{eq:n}) are
not satisfied.  We could try to overcome this difficulty by letting
$\delta = \delta(n) \to 0$ and $k = k(n) \to + \infty $ as $n $
increases. But doing this, we loose the exponential character of the
upper bound.  This could be considered as an illustration of the
result in Finesso et al. (1996) who proved that in the simpler case of
estimating the order of a Markov chain, it is not possible to have a
universal estimator with exponential bounds for the probability of
overestimation.

\section{Some final comments and bibliographic remarks}

Chains with memory of variable length were introduced in the
information theory literature by Rissanen (1983) as a universal system
for data compression.  Originally called by Rissanen {\sl tree
  machine}, {\it tree source}, {\sl context models}, etc., this class
of models recently became popular in the statistics literature under
the name of {\sl Variable Length Markov Chains (VLMC)}, coined by
B\"uhlmann and Wyner (1999).

Rissanen (1983) not only introduced the notion of variable memory
models but he also introduced the algorithm Context to estimate the
probabilistic context tree.  From Rissanen (1983) to Galves et
al. (2008), passing by Ron et al. (1996) and B\"uhlmann and Wyner
(1999), several variants of the algorithm Context have been presented
in the literature. In all the variants the decision to prune a branch
is taken by considering a {\em gain} function.

Rissanen (1983), B\"uhlmann and Wyner (1999)  and Duarte et al. (2006)
all defined the gain function in terms of the log likelihood ratio function. 
Rissanen (1983) proved the weak consistency
of the algorithm Context in the case where the contexts have a bounded
length. B\"uhlmann and Wyner (1999) proved the weak consistency of the
algorithm also in the finite case without assuming a prior known bound
on the maximal length of the memory but using a bound allowed to grow
with the size of the sample.

A different gain function was introduced in Galves et al. (2008),
considering differences between successive empirical transition
probabilities and comparing them with a given threshold $\delta .$ An
interesting consequence of the use of this different gain function was
obtained by Collet et al. (2007).  They proved that in the case of a
binary alphabet and when taking $\delta $ within a suitable interval,
it is possible to recover the context tree in the bounded case out
from a noisy sample where each symbol can be flipped with small
probability independently of the others.

The case of unbounded probabilistic context trees as far as we know
was first considered by Ferrari and Wyner (2003) who also proved a
weak consistency result for the algorithm Context in this more general
setting. The unbounded case was also considered by Csisz{\'a}r and
Talata (2006) who introduced a different approach for the estimation
of the probabilistic context tree using the Bayesian Information
Criterion (BIC) as well as the Minimum Description Length Principle
(MDL). We refer the reader to this last paper for a nice description
of other approaches and results in this field, including the context
tree maximizing algorithm by Willems et al. (1995). We also refer the
reader to Garivier (2006a, b) for recent and elegant results on the
BIC and the Context Tree Weighting Method (CTW).  Garivier (2006c) is
a very good presentation of models having memory of variable length,
BIC, MDL, CTW and related issues in the framework of information
theory.

With exception of Weinberger et al. (1995), the issue of the rate of
convergence of the algorithm estimating the probabilistic context tree
was not addressed in the literature until recently.  Weinberger et
al. (1995) proved in the bounded case that the probability that the
estimated tree differs from the finite context tree is summable as a
function of the sample size. Assuming weaker hypotheses than Ferrari
and Wyner (2003), Duarte et al. (2006) proved in the unbounded case
that the probability of error decreases as the inverse of the sample
size.

Leonardi (2007) obtained an upper bound for the rate of convergence of
penalized likelihood context tree estimators.  It showed that the
estimated context tree truncated at any fixed height approximates the
real truncated tree at a rate that decreases faster than the inverse
of an exponential function of the penalizing term. The proof mixes the
approaches of Galves et al. (2008) and Csisz{\'a}r and Talata (2006).

Several interesting papers have recently addressed the question of
classification of proteins and DNA sequences using models with memory
of variable length, which in bio-informatics are often called
prediction suffix trees (PST).  Many of these papers have been written
from a bio-informatics point of view focusing on the development of
new tools rather than being concerned with mathematically rigorous
proofs.  The interested reader can find a starting point to this
literature for instance in the papers by Bejerano et al. (2001),
Bejerano and Yona (2001), Eskin et al. (2000), Leonardi (2006) and
Miele et al. (2005).  The same type of analysis has been used
successfully to classification tasks in other domains like musicology
(Lartillot et al. 2003), linguistics (Selding et al. 2001), etc.

This presentation did not intend to be exhaustive and the bibliography
in many cases only gives a few hints about possible starting points to
the literature.  However, we think we have presented the state of the
art concerning the rate of convergence of context tree estimators.

In the introduction we said that {\sl Rissanen's ingenious idea was to
  construct a stochastic model that generalizes the notion of
  relevant domain (in biology or linguistics) to any kind of symbolic
  strings}.  Actually, God only knows what Jorma had in mind when he
invented this class of models. The French poet Paul Eluard wrote a
book called {\sl Les fr\`eres voyants}. This was the name given in the
middle-age to people guiding blind persons. So maybe Rissanen acted as
a {\sl fr\`ere voyant} using his intuition to push mathematics and
statistics into a challenging new direction.


\begin{thebibliography}{}

\bibitem{Bejeranoetal:2001}  G. Bejerano, Y. Seldin, H. Margalit
  and N. Tishby, ``Markovian domain fingerprinting: statistical
  segmentation of protein sequences'', {\it Bioinformatics}, vol. {17}, pp. 927--934, 2001.

\bibitem{bejeranoyona:2001} G. Bejerano and G. Yona, ``Variations on
  probabilistic suffix trees: statistical modeling and prediction of
  protein families'', {\it Bioinformatics}, vol. 17, Number 1, pp. 23--43, 2001.
           
\bibitem{buhlmann1999} P. B{\"u}hlmann and A.J. Wyner, ``Variable
  length Markov chains'', {\it Ann. Statist.}, vol. 27,  pp. 480--513, 1999.

\bibitem{colletetal:2007} P. Collet, A. Galves and F. Leonardi,
  ``Random perturbations of stochastic chains with unbounded variable
  length memory'', manuscript, can be downloaded from ArXiv:
  math/0707.2796., 2007.

\bibitem{comets:2002} F. Comets, R. Fern\'andez and P. Ferrari,
  ``Processes with long memory: Regenerative construction and perfect
  simulation'', {\it Ann. of Appl. Probab.}, vol. {12}, Number 3, pp. 921--943,
  2002.

\bibitem{csi:2002} I. Csisz{\'a}r, ``Large-scale typicality of {M}arkov
  sample paths and consistency of {MDL} order estimators. Special
  issue on Shannon theory: perspective, trends, and applications'', {\it IEEE
  Trans. Inform. Theory}, vol. {48}, Number 6, pp. 1616--1628, 2002.

\bibitem{csital:2006} I. Csisz{\'a}r and Z. Talata, ``Context tree
  estimation for not necessarily finite memory processes, via {BIC}
  and {MDL}'', {\it IEEE Trans. Inform. Theory}, vol. {52}, Number 3, pp. 1007--1016,
  2006.

\bibitem{dedeckerdoukhan:2003} J. Dedecker and P. Doukhan, ``A new
  covariance inequality and applications'', {\it Stochastic Process. Appl.},
  vol. {106}, Number 1, pp. 63--80, 2003.
  
\bibitem{dedecker:2005} J. Dedecker and C. Prieur, ``New dependence
  coefficients. Examples and applications to statistics'', {\it Probab.
  Theory Related Fields}, vol. {132}, pp. 203--236, 2005.

\bibitem{dgg:2006} D. Duarte, A. Galves and N.L. Garcia, ``Markov
  approximation and consistent estimation of unbounded probabilistic
  suffix trees'', {\it Bull. Braz. Math. Soc.}, vol. {37}, Number 4, pp. 581--592,
  2006.

\bibitem{es:2000 } E. Eskin, W.N. Grundy and Y. Singer, ``Protein
  family classification using sparse Markov transducers'', manuscript, 
  can be downloaded from http://citeseer.ist.psu.edu/328658.html,
  2000.
  
\bibitem{FG:2002} R. Fern{\'a}ndez and A. Galves, ``Markov
  approximations of chains of infinite order'', Fifth Brazilian School
  in Probability (Ubatuba, 2001), {\it Bull. Braz. Math. Soc. (N.S.)}, vol.{
    33}, Number 3, pp. 295--306, 2002.

\bibitem{fw:2003} F. Ferrari and A. Wyner, ``Estimation of general
  stationary processes by variable length {M}arkov chains'',
  {\it Scand. J. Statist.}, vol. {30}, Number 3, pp. 459--480, 2003.
  
\bibitem{finessoetal:1996} L. Finesso, C.-C. Liu P. and P. Narayan,
  ``The optimal error exponent for Markov order estimation'', {\it IEEE 
Trans. Inform. Theory}, vol. {42}, Number 5, pp. 1488--1497, 1996.

\bibitem{GL:2007} A. Galves and F. Leonardi, ``Exponential
  inequalities for empirical unbounded context trees'', manuscript, to
  appear in a
  special issue of {\it Progress in Probability}, V. Sidoravicius and M. E. Vares, Eds., 
  Birkh\"auser, 
  can be downloaded from ArXiv: math/0710.5900, 2007.

\bibitem{GMS:2008} A. Galves, V. Maume-Deschamps and B. Schmitt,
  ``Exponential inequalities for {VLMC} empirical trees'',
  {\it ESAIM. Probability and Statistics}, to appear 2008.

\bibitem{garivier:2006a} A. Garivier, ``Consistency of the unlimited
  BIC Context Tree estimator'', {\it IEEE Trans. Inform. Theory}, vol. {52},
  Number 10, pp. 4630--4635, 2006.

\bibitem{garivier:2006b} A. Garivier, ``Redundancy of the context-tree
  weighting method on renewal and Markov renewal processes'', {\it IEEE
  Trans. Inform. Theory}, vol. {52}, Number 12, pp. 5579--5586, 2006.

\bibitem{garivier:2006c} A. Garivier, ``Mod\`eles contextuels et
  alphabets infinis en th\'eorie de l'information'', PhD thesis,
  Universit\'e de Paris Sud, can be downloaded from
  http://www.math.u-psud.fr/~garivier/, 
  Nov. 2006.
 
\bibitem{haya:1977} T. Hayakawa, ``The likelihood ratio criterion and
  the asymptotic expansion of its distribution'',
  {\it Ann. Inst. Statist. Math.}, vol. {29}, Number 3, pp. 359--378, 1977.
              
\bibitem{lartillotetal:2003} O. Lartillot, S. Dubnov, G. Assayag and
  G. Bejerano, ``A system for computer music generation by learning and
  improvisation in a particular style'', {\it IEEE Computer J.}, vol. {36}, 
  Number 10, pp. 73--80, 2003.
              
\bibitem{leonardi:2006} F. Leonardi, ``A generalization of the PST
  algorithm: modeling the sparse nature of protein sequences'',
  {\it Bioinformatics}, vol. {22}, Number 7, pp. 1302--1307, 2006.

\bibitem{leonardi:2007} F. Leonardi, ``Rate of convergence of penalized
  likelihood context tree estimators'', manuscript,  can be downloaded from ArXiv:
  math/0701810v2, 2007.

\bibitem{mar:shi:1994} K. Marton and P.C. Shields, ``Entropy and the
  consistent estimation of joint distributions'', {\it Ann. Probab.}, vol. {22},
  Number 2, pp. 960--977, 1994.

\bibitem{mar:shi:1996} K. Marton and P.C. Shields, ``Correction:
  ``{E}ntropy and the consistent estimation of joint distributions''
  [{A}nn.\ {P}robab.\ {\bf 22} (1994), no.\ 2, 960--977; {MR}1288138
  (95g:94004)]'', {\it Ann. Probab.}, vol. {24}, Number 1, pp. 541--545, 1996.

\bibitem{maume:2006} V. Maume-Deschamps, {\it Exponential
    inequalities and estimation of conditional probabilities}, Lecture
  Notes in Statist., vol. 187, Heidelberg: Springer, 2006.

\bibitem{miel:2005} V. Miele, P. Y. Bourguignon, D. Robelin, G. Nuel
  and H. Richard, ``seq++: a package for biological sequences analysis
  with a range of Markov-related models'', {\it BioInformatics}, vol. {21},
  Number 11, pp. 2783--2874, 2005.

\bibitem{rissanen1983} J. Rissanen, ``A universal data compression
  system'', {\it IEEE Trans. Inform. Theory}, vol. {29}, Number 5, pp. 656--664,
  1983.

\bibitem{ronetal:1996} D. Ron, Y. Singer and N. Tishby, ``The power
  of amnesia: learning probabilistic automata with variable memory
  length'', {\it Machine Learning}, vol. {25}, pp. 117--149, 1996.

\bibitem{seldingetal:2001} Y. Seldin, G. Bejerano and N. Tishby,
  ``Unsupervised sequence segmentation by a mixture of variable memory
  Markov models'', manuscript, can be downloaded from
  http://citeseer.ist.psu.edu/449505.html, 2001.

\bibitem{vaart:1998} A.W. van der Vaart, {\it Asymptotic
    statistics}, Cambridge Series in Statistical and Probabilistic
  Mathematics, Cambridge: Cambridge University Press, 1998.

\bibitem{wrf:1995} M.J. Weinberger, J. Rissanen and M. Feder, ``A
  universal finite memory source'', {\it IEEE Trans. Inform. Theory}, vol. {41},
  Number 3, pp. 643--652, 1995.

\bibitem{willemsetal} E. M. J. Willems, Y. M. Shtarkov and T. J. Tjalkens, 
  ``The Context-Tree weighting method: Basic properties'', {\it IEEE
  Trans. Inform. Theory}, vol. {41}, pp. 653--664, 1995.

\end{thebibliography}

\vskip30pt

Antonio Galves

Instituto de Matem\'atica e Estat\'{\i}stica

Universidade de S\~ao Paulo

BP 66281

05315-970 S\~ao Paulo, Brasil

e-mail: {\tt galves@ime.usp.br}

\bigskip

Eva L\"ocherbach

Universit\'e Paris-Est

LAMA -- UMR CNRS 8050

61, Avenue du G\'en\'eral de Gaulle 

94000 Cr\'eteil, France

e-mail: {\tt locherbach@univ-paris12.fr}
\end{document}